\documentclass[CJK,11pt,twoside, openany]{ctexart}
\usepackage{latexsym}
\usepackage{amsthm}
\usepackage{bbm}
 \usepackage{mathrsfs}
 \usepackage{fancyhdr}
\usepackage{mathrsfs}
\usepackage{cancel}
\usepackage{color}
\usepackage{multirow}
\usepackage{eso-pic}
\usepackage{dsfont}
\usepackage{amssymb}
\def\nnb{\nonumber}
\def\ds{\displaystyle}
\def\cd{\cdot}

\def\all{  \, \forall \, }

\newcommand{\refeq}[1]{~$(\ref{#1})$}
\newcommand{\eqref}[1]{~$(\ref{#1})$ }

\def\endpf{\hfill$\Box$\vspace{0.4cm}}
\def\eqif{\, {\rm if}\, \,}

\def\eqin{ \, {\rm in } \, \,}
\def\eqae{ \, {\rm a.e. } \,\, }

\newcommand{\wi}[1]{\mbox{ weakly  in}\, #1}

\def\ol{\overline}
\def\Ga{\alpha}
\def\Gb{\beta}

\def\Gvp{\varphi}

\def\bY{\ol{Y}}

\def\bu{{\bar u}}

\def\by{{\bar y}}

\def\mcP{{\mathscr P}}

\def\mcU{{\mathscr U}}

\def\tiQ{\widetilde  Q}

\def\qq{\qquad}
\def\q{\quad}
\newcommand{\doo}[2]{{d {#1} \over d {#2}}}

\newcommand{\pri}{{\prime}}
\newcommand{\prii}{{\prime\prime}}

\def\IR{\mathds{R}}
\def\defeq{\buildrel \triangle \over =}

\newcommand{\set}[1]{\left\{#1\right\}}

\def\thebibliography#1{\center{\bf\normalsize References}\list
 {[\arabic{enumi}]}{\settowidth\labelwidth{[#1]}\leftmargin\labelwidth
 \advance\leftmargin\labelsep
 \usecounter{enumi}}
 \def\newblock{\hskip .11em plus .33em minus .07em}
 \sloppy\clubpenalty4000\widowpenalty4000
 \sfcode`\.=1000\relax}

\makeatletter
\def\cleardoublepage{\clearpage\if@twoside \ifodd\c@page\else
   \hbox{}\thispagestyle{empty}\newpage\addtocounter{page}{-1}
   \if@twocolumn\hbox{}\newpage\fi\fi\fi}
\makeatother
\usepackage{graphicx}
\newtheorem{Definition}{Definition}[section]
\newtheorem{Theorem}[Definition]{Theorem}
\newtheorem{Lemma}[Definition]{Lemma}
\newtheorem{Corollary}[Definition]{Corollary}
\newtheorem{Proposition}[Definition]{Proposition}
\newtheorem{Remark}{Remark}[section]

\begin{document}
\title{\bf Sharp Constants for Inequalities of Poincar\'e Type: An Application of Optimal Control Theory \thanks{This work was supported in part by 973 Program (No. 2011CB808002)  and  NSFC (No. 61074047).}}

\author{Hongwei Lou\footnote{School of Mathematical Sciences, and LMNS, Fudan University, Shanghai
200433, China (Email: \texttt{hwlou@fudan.edu.cn}).}
}

\date{}

\maketitle

\begin{quote}
\footnotesize {\bf Abstract.} Sharp constants for an inequality of Poincar\'e type is studied.
The problem is solved by using optimal control theory.

\textbf{Key words and phrases.} Sharp constant,  inequality of Poincar\'e type, optimal control.

\textbf{AMS subject classifications.} 26D10, 46E35, 49K15.
\end{quote}

\normalsize

\def\theequation{1.\arabic{equation}}
\setcounter{equation}{0} 
\setcounter{Definition}{0} \setcounter{Remark}{0}
 \vspace{6mm}

\section{Introduction.}
Denote by $W^{1,2}(-1,1)$ the Sobolev space of all real-valued functions $f(\cd)$ that are absolutely continuous
on the closed interval $[-1, 1]$ and such that $f^\pri(\cd)\in L^2(-1,1)$.
Let $m\geq 1$ be an integer. Denote by $W(1,2,m)$ the space
\begin{equation}\label{E101}
W(1;2,m)\defeq \set{y(\cd)\in W^{1,2}(-1,1)\Big|  \int^1_{-1}t^ky(t)\, dt=0, \q (0\leq k\leq m-1)}.
\end{equation}
G. A. Kalyabin considered in \cite{Kal} the following problem:

\textbf{Problem ($B_x$):} Fix $x\in [-1,1]$, find the best constant $B_m(x)$ such that the following inequality holds:
\begin{equation}\label{E102}
   |y(x)|\leq B_m(x)\Big(\int^1_{-1}|y^\pri(t)|^2\, dt\Big)^{1\over 2},\qq \all y(\cd)\in W(1,2,m).
\end{equation}
It is proved in \cite{Kal} that
\begin{equation}\label{E103}
B_m^2(\pm 1)={2\over m(m+2)}, \q B_1^2(x)={1+3x^2\over 6}, \q B_2^2(0)={1\over 6}
\end{equation}
and the extremal functions for the case $x=1$ is
\begin{equation}\label{E104}
C\Big((m+2)P_m(t)+mP_{m+1}(t)\Big),
\end{equation}
 where $C$ is a constant and
\begin{equation}\label{E105}
 P_k(t)={1\over 2^k k!}{d^k\over dt^k}(t^2-1)^k, \qq k=0,1,2,\ldots
\end{equation}
are the classical Legendre polynomials. For notation simplicity, we denote
\begin{equation}\label{E106}
 p_k(t)={1\over 2^k k!}{d^{k-2}\over dt^{k-2}}(t^2-1)^k, \qq k=2,3,\ldots
\end{equation}
and
\begin{equation}\label{E107}
\Ga_k(x)={\int^1_{-1}P_k(s)\, |s-x|\, ds\over \int^1_{-1}|P_k(s)|^2\, ds}=
\left\{\begin{array}{ll}{x^2+1\over 2}, & \eqif k=0, \\
 {x^3-3x\over 2}, & \eqif k=1, \\
(2k+1)p_k(x), &  \eqif k\geq 2.
\end{array}\right.
\end{equation}

In this paper, we will solve Problem ($B_x$) completely by the help of optimal control theory.
Since the cases $x=\pm 1$ were solved in \cite{Kal}, we mainly consider cases of $x\in (-1,1)$. We have
\begin{Theorem}\label{T101}
Assume $m\geq 1$ and  $x\in (-1,1)$.
Then,  $y(\cd)\in W(1,2,m)$ is an extremal function to
 Problem $(B_x)$ if and only if
\begin{equation}\label{E108}
y(t)=C\Big[c(x)\Big(Q_{m+1}(t)-|t-x|\Big)+1\Big],\qq t\in [-1,1],
\end{equation}
where $C$ is a constant,
\begin{equation}\label{E109}
Q_{m+1}(t)=a(x)+\sum^{m-1}_{k=1}\Ga_k(x)\, P_k(t)+\Ga(x) P_m(t)+\Gb(x) P_{m+1}(t), \q t\in [-1,1],
\end{equation}
\begin{equation}\label{E110}
c(x)={2\over x^2+1-2a(x)}
\end{equation}
and
$a(x),  \Ga(x), \Gb(x)$ are characterized by
\begin{equation}\label{E111}
Q_{m+1}(x)=0, \q Q^\pri_{m+1}(-1)=-1, \q Q^\pri_{m+1}(1)=1.
\end{equation}
While the sharp constant of inequality \refeq{E102} is
\begin{equation}\label{E112}
B_m(x)= {1\over \sqrt{2c(x)}}, \qq x\in (-1,1).
\end{equation}
\end{Theorem}
More precisely, we have
\begin{Corollary}\label{T102}
Assume $m=1$ and $x\in (-1,1)$.
Then
\begin{equation}\label{E113}
B_1(x)=\sqrt{3x^2+1\over 6}
\end{equation}
and $y(\cd)\in W(1,2,1)$ is an extremal function to
 Problem $(B_x)$ if and only if
 \begin{equation}\label{E114}
y(t)=C\Big[{3\over 3x^2+1}\Big({t^2-x^2\over 2}-|t-x|\Big)+1\Big], \q t\in [-1,1].
\end{equation}
\end{Corollary}
\begin{Corollary}\label{T103}
Assume $m=2$ and $x\in (-1,1)$.
Then
\begin{equation}\label{E115}
B_2(x)=\sqrt{8-21x^2+30x^4-5x^6\over 48}
\end{equation}
and $y(\cd)\in W(1,2,2)$ is an extremal function to
 Problem $(B_x)$ if and only if
 \begin{equation}\label{E116}
y(t)=C\Big[{4+33x^2-30x^4+5x^6+8P_2(t)-2(x^3-3x)(P_3(t)-6P_1(t))-24|t-x|\over 8-21x^2+30x^4-5x^6}+1\Big].
\end{equation}
\end{Corollary}
\begin{Corollary}\label{T104}
Assume $m=2n+1$, $n\geq 1$ and $x\in (-1,1)$.
Then
\begin{eqnarray}\label{E117}
\nnb && B_{2n+1}(x) =\Big\{ {(x^2-1)^2\over 4}+{1\over 2}\sum^{2n}_{k=2} (2k+1) p_k(x)\, P_k(x)\\
\nnb  && -{1\over 2(n+1)(2n+1)}\Big[ {x^3-3x \over 2}+\sum^{n-1}_{k=1} (k+1)(2k+1)(4k+3) p_{2k+1}(x)\Big]P_{2n+1}(x)\\
 && +{1\over 2(n+1)(2n+3)}\Big[1-\sum^n_{k=1} k(2k+1)(4k+1) p_{2k}(x)\Big]P_{2n+2}(x) \Big\}^{1\over 2}.
\end{eqnarray}
In particular,
 \begin{equation}\label{E118}
B_3(x)={\sqrt {297+1260x^2-5370x^4+5900x^6-1575x^8}\over 16\sqrt{15}},
\end{equation}
 \begin{equation}\label{E119}
B_5(x)={\sqrt {1375+8400x^2-95025x^4-357560x^6-597555x^8+448056x^{10}-121275x^{12}}\over 16\sqrt{105}}.
\end{equation}
\end{Corollary}
\begin{Corollary}\label{T105}
Assume  $m=2n+2$, $n\geq 1$ and $x\in (-1,1)$.
Then
\begin{eqnarray}\label{E120}
\nnb && B_{2n+2}(x) =\Big\{ {(x^2-1)^2\over 4}+{1\over 2}\sum^{2n+1}_{k=2} (2k+1) p_k(x)\, P_k(x)\\
\nnb &&+ {1\over 2(n+1)(2n+3)}\Big[1-\sum^n_{k=1} k(2k+1)(4k+1) p_{2k}(x)\Big] P_{2n+2}(x)\\
&& -{1\over 2(n+2)(2n+3)}\Big[ {x^3-3x \over 2}+\sum^n_{k=1} (k+1)(2k+1)(4k+3) p_{2k+1}(x)\Big]P_{2n+3}(x) \Big\}^{1\over 2}.
\end{eqnarray}
In particular,
 \begin{equation}\label{E121}
B_4(x)={\sqrt {297-1440x^2+9030x^4-20860x^6+18585x^8-5292x^{10}}\over 16\sqrt{15}},
\end{equation}
 \begin{eqnarray}\label{E122}
\nnb B_6(x)&=& {1\over 32\sqrt{42}}\Big(2200-15225x^2+211050x^4-1162455x^6+3017700x^8\\
&& -3977127x^{10}+2562714x^{12}-637065x^{14}\Big)^{1\over 2}.
\end{eqnarray}
\end{Corollary}
When $x=0$, we have
\begin{Corollary}\label{T106} It holds that
 \begin{equation}\label{E123}
B_1(0)=B_2(0)={\sqrt 6\over 6},\q B_3(0)=B_4(0)={3\sqrt{55}\over 80}
\end{equation}
and
\begin{eqnarray}\label{E124}
\nnb && B_{2n+1}(0)=B_{2n+2}(0) \\
\nnb &=& \Big\{ {3\over 32} -{1\over 2}\sum^n_{k=2} (4k+1) {(2k-3)!!\over (2k+2)!!}\cd{(2k-1)!!\over (2k)!!}+{7\cd (-1)^n\over 16(n+1)(2n+3)}{(2n+1)!!\over (2n+2)!!}\\
 &&-\sum^n_{k=2}{(-1)^{n-k}k(2k+1)(4k+1)\over 2(n+1)(2n+3)}{(2k-3)!!\over (2k+2)!!}\cd{(2n+1)!!\over (2n+2)!!}\Big\}^{1\over 2},\q n\geq 2.
\end{eqnarray}
In particular,
 \begin{equation}\label{E125}
B_6(0)=B_6(0)=\sqrt{ 275\over 5376},\q B_7(0)=B_8(0)=\sqrt{45325\over 1179648}.
\end{equation}
\end{Corollary}

\bigskip

\def\theequation{2.\arabic{equation}}
\setcounter{equation}{0} 
\setcounter{Definition}{0} \setcounter{Remark}{0}\section{Transmit Problem ($B_x$) to Optimal Control Problem.}

We introduce the equivalent optimal control problem to Problem $(B_x)$.
Let $\mcU=L^2(-1,1)$. We set the following control system:
\begin{equation}\label{E201}
 \doo{} t  \pmatrix{y(t)\cr w(t)\cr z_0(t) \cr z_1(t)\cr \vdots\cr z_{m-1}(t)}=
    \pmatrix{u(t)\cr u(t)\chi_{(-1,x)}(t)\cr y(t) \cr ty(t)\cr \vdots \cr t^{m-1}y(t)}, \qq t\in (-1,1)
 \end{equation}
and the state constraints
\begin{equation}\label{E202}
y(-1)=w(-1), \q w(1)=1, \q z_k(\pm 1)=0, \q (k=0,1,\ldots,m-1).
 \end{equation}
Let
\begin{equation}\label{E203}
\mcP_{ad}=\set{(Y(\cd),u(\cd))\in \Big(W^{1,2}(-1,1)\Big)^{m+2}\times \mcU\Big|(Y(\cd),u(\cd))\q \mbox{satisfies \refeq{E201}---\refeq{E202}}}
\end{equation}
and
\begin{equation}\label{E204}
\mcU_{ad} =\set{u(\cd)|(Y(\cd),u(\cd))\in \mcP_{ad}},
\end{equation}
 where
\begin{equation}\label{E205}
Y(\cd)=\pmatrix{y(\cd)\cr w(\cd)\cr z_0(\cd) \cr z_1(\cd)\cr \vdots\cr z_{m-1}(\cd)}.
\end{equation}
Our optimal control problem corresponding to Problem $(B_x)$ is

\textbf{Problem ($C_x$).} Let $x\in (-1,1)$. Find $(\bY(\cd),\bu(\cd))\in \mcP_{ad}$  such that
\begin{equation}\label{E206}
\int^1_{-1}\bu^2(t)\, dt=\inf_{(Y(\cd),u(\cd))\in \mcP_{ad}}\int^1_{-1}u^2(t)\, dt.
\end{equation}

It is obvious that $(Y(\cd),u(\cd))\mapsto y(\cd)$ is a bijection from $\mcP_{ad}$ to $\set{f(\cd)\in W(1,2,m)|f(x)=1}$. Then one can easily see that
\begin{equation}\label{E207}
B_m(x)=\Big(\inf_{(Y(\cd),u(\cd))\in \mcP_{ad}}\int^1_{-1}u^2(t)\, dt\Big)^{-{1\over 2}}.
\end{equation}
Therefore, we can solve Problem ($B_x$) by solving Problem ($C_x$).

\bigskip

\def\theequation{3.\arabic{equation}}
\setcounter{equation}{0} 
\setcounter{Definition}{0} \setcounter{Remark}{0}\section{Proof of Theorem \ref{T101}.}

We give the following lemma first.
\begin{Lemma}\label{T301} Let $n\geq 1$, $c\in \IR$,  $x\in (-1,1)$, $Q(\cd)$ is an $(n+1)$-th degree polynomial satisfying
\begin{equation}\label{E301}
Q(x)=0, \q Q^\pri(1)=1, \q Q^\pri(-1)=-1
\end{equation}
and
\begin{equation}\label{E302}
\int^1_{-1} t^k\Big[c\Big(Q(t)-|t-x|\Big)+1\Big]\, dt=0, \qq \all k=0,1,\ldots,n-1.
\end{equation}
Then
\begin{equation}\label{E303}
\int^1_{-1}  \Big[{\partial \over \partial t} \Big(Q(t)-|t-x|\Big)\Big]^2\, dt={2\over c}.
\end{equation}
\end{Lemma}
\proof Noting that $Q^\prii(\cd)$ is an $(n-1)$-th degree polynomial, by \refeq{E302}, we have
\begin{equation}\label{E304}
\int^1_{-1} Q^\prii(t)\Big[c\Big(Q(t)-|t-x|\Big)+1\Big]\, dt=0.
\end{equation}
Therefore,
\begin{eqnarray}\label{E305}
\nnb && \int^1_{-1} \Big[{\partial \over \partial t} \Big(Q(t)-|t-x|\Big)\Big]^2\, dt\\
\nnb &=&\int^x_{-1}  \Big(Q^\pri(t)+1\Big)^2\, dt+\int^1_x  \Big(Q^\pri(t)-1\Big)^2\, dt\\
\nnb &=&\int^1_{-1}  \Big(Q^\pri(t)\Big)^2\, dt+2+2\int^x_{-1}  Q^\pri(t)\, dt-2\int^1_x Q^\pri(t)\, dt\\
\nnb &=& Q(1)Q^\pri(1)-Q(-1)Q^\pri(-1)-\int^1_{-1}   Q^\prii(t)Q(t)\, dt+2+4Q(x)-2Q(-1)-2Q(1)\\
\nnb &=&  \int^1_{-1} Q^\prii(t)\Big({1\over c}-|t-x|\Big)\, dt+2- Q(-1)- Q(1)\\
\nnb &=& {2\over c}+ \int^x_{-1} Q^\prii(t) (t-x)\, dt+\int^1_x Q^\prii(t) (x-t)\, dt+2- Q(-1)- Q(1)\\
\nnb &=& {2\over c}-(1+x)- \int^x_{-1} Q^\pri(t)\, dt-(1-x)+\int^1_x Q^\pri(t)\, dt+2- Q(-1)- Q(1)\\
&=& {2\over c}.
\end{eqnarray}
\endpf

Now, we list some properties of Legendre polynomials.
First, we mention some formulae that are useful to us. We can get easily that
\begin{equation}\label{E306}
P_k(0)
=\left\{\begin{array}{ll} 0, & \eqif k=2n+1, \, n\geq 0, \\
1,  & \eqif k=0, \\
(-1)^n{(2n-1)!!\over (2n)!!},  & \eqif k=2n, \, n\geq 1,
\end{array}\right.
\end{equation}
\begin{equation}\label{E307}
p_k(0)
=\left\{\begin{array}{ll} 0, & \eqif k=2n+1, \, n\geq 1, \\
{1\over 8},  & \eqif k=2, \\
(-1)^{n-1}{(2n-3)!!\over (2n+2)!!},  & \eqif k=2n, \, n\geq 2,
\end{array}\right.
\end{equation}

\begin{equation}\label{E308}
P_k(-1)=(-1)^k, \q P_k(1)=1,
\end{equation}
\begin{equation}\label{E309}
P^\pri_k(-1)=(-1)^{k-1}{k(k+1)\over 2}, \q P^\pri_k(1)={k(k+1)\over 2},
\end{equation}
\begin{eqnarray}\label{E310}
&&\int^1_{-1}|P_k(t)|^2\, dt = {2\over 2k+1}.
\end{eqnarray}
and
\begin{equation}\label{E311}
\int^1_{-1}(P^\pri_k(t))^2\, dt=P^\pri_k(1)P_k(1)-P^\pri_k(-1)P_k(-1)=k(k+1).
\end{equation}
\begin{equation}\label{E312}
\int^1_{-1}P^\pri_k(t)P^\pri_{k+1}(t)\, dt=P^\pri_k(1)P_{k+1}(1)-P^\pri_k(-1)P_{k+1}(-1)=0.
\end{equation}
\begin{equation}\label{E313}
\int^1_{-1} P_k(t)|t-x|\, dt=\left\{\begin{array}{ll}x^2+1, & \eqif k=0, \\
 {1\over 3}x^3-x, & \eqif k=1, \\
2 p_k(x), &  \eqif k\geq 2.
\end{array}\right.
\end{equation}

\bigskip

We turn to prove Theorem \ref{T101}.

\textbf{Proof of Theorem \ref{T101}.}

\textbf{I. Existence of optimal pair.} One can prove directly that the sharp constant $B_m(x)$ is attainable, i.e., there is a nontrivial $\by(\cd)\in W(1,2,m)$ such that \begin{equation}\label{E314}
|\by(x)|=B_m(x)\Big(\int^1_{-1} |y^\pri(t)|^2\, dt\Big)^{-{1\over 2}}.
\end{equation}
Now, we give an optimal control version of this fact.

Let $(Y_j(\cd),u_j(\cd))\in \mcP_{ad}$ be a minimizing sequence of Problem ($C_x$). That is
\begin{equation}\label{E315}
\lim_{j\to +\infty}\int^1_{-1} u_j^2(t)\, dt=\inf_{(Y(\cd),u(\cd))\in \mcP_{ad}}\int^1_{-1}u^2(t)\, dt.
\end{equation}
Then $\ds u_j(\cd)$ is bounded in $L^2(-1,1)$. By Eberlein-Shmulyan Theorem (see \cite{Yo}, for example), we can suppose that
\begin{equation}\label{E316}
u_j(\cd)\to \bu(\cd), \qq \wi L^2(-1,1)
\end{equation}
for some $\bu(\cd)\in L^2(-1,1)$. Then using the state constraints
\refeq{E202}, we can easily prove that
\begin{equation}\label{E317}
Y_j(\cd)\to \bY(\cd), \qq \eqin (C[-1,1])^{m+2}
\end{equation}
for some $\bY(\cd)\in (W^{1,2}(-1,1))^{m+2}$ and $(\bY(\cd),\bu(\cd))\in \mcP_{ad}$.  Moreover,
\begin{equation}\label{E318}
\int^1_{-1}\bu^2(t)\, dt\leq \lim_{j\to +\infty}\int^1_{-1} u_j^2(t)\, dt.
\end{equation}
Therefore $(\bY(\cd),\bu(\cd))$ is a solution to Problem ($C_x$). We call it as an optimal pair of  Problem ($C_x$).

\textbf{II. Pontryagin's maximum principle for the optimal pair.} We apply  Pontryagin's maximum principle. By optimal control theory, the optimal pair $(\bY(\cd),\bu(\cd))$ satisfies the following Pontryagin's maximum principle (see \cite{Ber} and Remark \ref{R302}): there exists a
$\Gvp^0\leq 0$ and a solution to the following conjugate equation
\begin{equation}\label{E319}
\ds \doo{} t  \pmatrix{\Gvp(t)\cr \zeta(t)\cr \psi_0(t) \cr \psi_1(t)\cr \vdots\cr \psi_{m-1}(t)}=
                  \pmatrix{\ds -\sum^{m-1}_{j=0}t^j\psi_j(t)\cr 0\cr 0 \cr 0\cr \vdots\cr 0},\qq t\in (-1,1)
\end{equation}
such that the following conditions hold:

(i) the non-trivial condition
\begin{equation}\label{E320}
(\Gvp^0,\Gvp(\cd), \zeta(\cd), \psi_0(\cd) , \ldots, \psi_{m-1}(\cd))\ne 0,
\end{equation}

(ii) the maximum condition
\begin{eqnarray}\label{E321}
\nnb && \Gvp^0 \bu^2(t)+\Big(\Gvp(t)+\zeta(t)\chi_{(-1,x)}(t)\Big)\bu(t)\\
&=& \max_{u\in \IR}\Big[\Gvp^0 u^2+\Big(\Gvp(t)+\zeta(t)\chi_{(-1,x)}(t)\Big) u\Big],\qq \eqae t\in (-1,1),
\end{eqnarray}

(iii) the conversation condition
\begin{equation}\label{E322}
\Gvp(-1)+\zeta(-1)=0, \q \Gvp(1)=0.
\end{equation}

\textbf{III. Analyze.}  By \refeq{E319}, $\zeta(\cd)\equiv\zeta, \psi_0(\cd)\equiv\psi_0, \ldots,\psi_{m-1}(\cd)\equiv\psi_{m-1}$ are constants\footnote{here and hereafter, constants may depend on $x$.} and
\begin{equation}\label{E323}
\Gvp^\pri(t)=-\sum^{m-1}_{j=0}\psi_j t^j, \qq t\in (-1,1).
\end{equation}
Then $\Gvp(\cd)$ is an $m$-th degree polynomial.

By \refeq{E321}, we have
\begin{equation}\label{E324}
2\Gvp^0\bu(t)+\Gvp(t)+\zeta\chi_{(-1,x)}=0, \qq\eqae t\in (-1,1).
\end{equation}
If $\Gvp^0=0$, then since $x\in (-1,1)$, we get
\begin{equation}\label{E325}
\Gvp(\cd)\equiv 0,\q \zeta=0.
\end{equation}
This contradicts to the non-trivial condition \refeq{E320}. Therefore, we must have $\Gvp^0<0$. Without loss of generality, we can suppose that
$\Gvp^0=-{1\over 2}$. Then it follows from \refeq{E324} that
\begin{equation}\label{E326}
\bu(t) =\Gvp(t)+\zeta\chi_{(-1,x)}(t), \qq t\in (-1,1).
\end{equation}
Combining with \refeq{E322}, we see that the corresponding function $\by(\cd)$ is deferential continuous on $[-1,x)\bigcup (x,1]$ and
\begin{equation}\label{E327}
\by^\pri(\pm 1) =\bu(\pm 1)=0.
\end{equation}
Moreover, $\by(\cd)$ can be expressed as
\begin{equation}\label{E328}
\by(t)=-c(x)|t-x|+\tiQ_{m+1}(t),\qq t\in [-1,1],
\end{equation}
where $c(x)={\zeta\over 2}$ and $\tiQ_{m+1}(\cd)$ is an $(m+1)$-th degree polynomial.

We claim that $c\ne 0$. Otherwise, $c=0$ and it follows from \refeq{E328} and
\begin{equation}\label{E329}
\int^1_{-1}t^k\by(t)\, dt=0, \qq k=0,1,2,\ldots,m-1
\end{equation}
that
\begin{equation}\label{E330}
\tiQ_{m+1}(t)=c_m P_m(t) +c_{m+1} P_{m+1}(t), \qq k=0,1,2,\ldots,m-1
\end{equation}
for some constant $c_m,c_{m+1}$.

Then,  \refeq{E327} and  \refeq{E309} imply $c_m=c_{m+1}=0$. This contradicts to the nontrivial condition.
Therefore $c\ne 0$ and we can rewrite $\by(\cd)$ as
\begin{equation}\label{E331}
\by(t)=c(x) \Big( Q_{m+1}(t)-|t-x|\Big)+1,\qq t\in [-1,1],
\end{equation}
where $Q_{m+1}(\cd)$ is an $(m+1)$-th degree polynomial such that
\begin{equation}\label{E332}
Q_{m+1}(x)=0, \q Q^\pri_{m+1}(-1)=-1, \q Q^\pri_{m+1}(1)=1.
\end{equation}

\textbf{IV. Conclusion.}
By  \refeq{E329}, we can get that
\begin{equation}\label{E333}
Q_{m+1}(t)=a(x)+\sum^{m-1}_{k=1}\Ga_k(x)\, P_k(t)+\Ga(x) P_m(t)+\Gb(x) P_{m+1}(t),
\end{equation}
where $\Ga_k(x)$ is defined by \refeq{E107}. Moreover, we can determine $a(x),\Ga(x)$ and $\Gb(x)$ by \refeq{E332}. Finally, using \refeq{E329} again, we get that
\begin{equation}\label{E334}
c(x)={2\over x^2+1-2a(x)}.
\end{equation}
Then Theorem \ref{T101} follows from \refeq{E207} and Lemma \ref{T301}.\endpf

\begin{Remark}\label{R301} If $x=-1$, instead of \refeq{E331}---\refeq{E332}, we could get that
\begin{equation}\label{E335}
\by(t)=Q_{m+1}(t),\qq t\in [-1,1]
\end{equation}
with
\begin{equation}\label{E336}
Q_{m+1}(-1)=1, \q Q^\pri_{m+1}(1)=0
\end{equation}
and
\begin{equation}\label{E337}
Q_{m+1}(t)=\Ga P_m(t)+\Gb P_{m+1}(t),\qq t\in [-1,1].
\end{equation}
The above equations imply the results got in \cite{Kal} for $x=\pm 1$.
\end{Remark}
\begin{Remark}\label{R302} \if{We should indicate that it seems not easy to
find a reference containing a Pontryagin maximum principle that can be used directly here.
For technical reason, to yield a  Pontryagin maximum principle, people usually assume that the right hand term of the system equation $($Equation \refeq{E201} in this paper$)$ is controlled uniformly by an integrable function. However, the right hand term of \refeq{E201} can not be controlled uniformly by an integrable function for all $u(\cd)\in \mcU_{ad}$.

Nevertheless, we are sure that researchers in optimal control theory know that
the Pontryagin maximum principle we used in part \textbf{II} is valid.

On the other hand, if we replace $\mcU_{ad}$ by
$$
\widetilde{\mcU}_{ad}\equiv \set{u(\cd)\in \mcU_{ad}\Big| \|u(\cd)-\bu(\cd)\|_{L^\infty(-1,1)}},
$$
then Pontryagin maximum principles in many references can be used directly here and the discussion in the above proof can go still smoothly.
}\fi
Strictly speaking, mainly because there is not a $\xi(\cd)\in L^1(-1,1)$ which can control $u(\cd)\in \mcU_{ad}$ uniformly, we can not use
the Pontryagin maximum principles given in \cite{Ber} directly.
Nevertheless,  the Pontryagin maximum principle we used in part \textbf{II} is valid.

If we want use the results in \cite{Ber} (or Pontryagin maximum principles in many other references) directly, we need only
replace $\mcU_{ad}$ by
$$
\widetilde{\mcU}_{ad}\equiv \set{u(\cd)\in \mcU_{ad}\Big| \|u(\cd)-\bu(\cd)\|_{L^\infty(-1,1)}<1}.
$$
And the the proof of Theorem \ref{T101} could be same to what we have just given.
\end{Remark}

\bigskip

\def\theequation{4.\arabic{equation}}
\setcounter{equation}{0} 
\setcounter{Definition}{0} \setcounter{Remark}{0}\section{Results for some spacial cases.}

We prove Corollaries \ref{T102}---\ref{T106} in this section.

\textbf{Proof of Corollary \ref{T102}.}

By \refeq{E107} and \refeq{E109}, we have
\begin{equation}\label{E401}
Q_2(t)= {1\over 3}\Big(P_2(t)-P_2(x)\Big)={t^2-x^2\over 2}.
\end{equation}
That is
\begin{equation}\label{E402}
a(x)=-{1\over 3} P_2(x)={1-3x^2\over 6}.
\end{equation}
Then
\begin{equation}\label{E403}
c(x)\equiv {2\over x^2+1-2a(x)}={3\over 3x^2+1}.
\end{equation}
Therefore, the extremal function to Problem ($B_x$) is $C\by(\cd)$ with
\begin{equation}\label{E404}
\by(x)={3\over 3x^2+1}\Big({t^2-x^2\over 2}-|t-x|\Big)+1.
\end{equation}
While
\begin{equation}\label{E405}
B_1(x)= {1\over \sqrt{2c(x)}}=\sqrt{3x^2+1\over 6}.
\end{equation}
\endpf

\textbf{Proof of Corollary \ref{T103}.}

By \refeq{E107} and \refeq{E109}, we have
\begin{equation}\label{E406}
 Q_3(t) =   a(x)+{x^3-3x \over 2}\, P_1(t)+\Ga(x) P_2(t)+\Gb(x) P_3(t).
\end{equation}
Then it follows easily from  \refeq{E111} that
\begin{equation}\label{E407}
\Ga(x)  =  {1\over 3}, \q \Gb(x)=  -{x^3-3x \over 12}, \q a(x)={4+33x^2-30x^4+5x^6\over 24}.
\end{equation}
Thus
\begin{equation}\label{E408}
c(x)  = {2\over x^2+1-2a}={24\over 8-21x^2+30x^4-5x^6}.
\end{equation}
Therefore
\begin{equation}\label{E409}
B_2(x)= {1\over \sqrt{2c(x)}}=\sqrt{8-21x^2+30x^4-5x^6\over 48}
\end{equation}
and the  extremal functions to Problem ($B_x$) are
\begin{equation}\label{E410}
C\Big({4+33x^2-30x^4+5x^6+8P_2(t)-2(x^3-3x)(P_3(t)-6P_1(t))-24|t-x|\over 8-21x^2+30x^4-5x^6}+1\Big).
\end{equation}
\endpf

\textbf{Proof of Corollary \ref{T104}.} 

By \refeq{E107} and \refeq{E109}, we have
\begin{eqnarray}\label{E411}
\nnb Q_{2n+2}(t) &=& a(x)+ {x^3-3x \over 2}P_1(t)+\sum^{2n}_{k=2} (2k+1) p_k(x)\, P_k(t)\\
&& +\Ga(x) P_{2n+1}(t)+\Gb(x) P_{2n+2}(t).
\end{eqnarray}
Then by  \refeq{E111},
\begin{eqnarray}
\label{E412} \Ga(x) &=& -{1\over (n+1)(2n+1)}\Big[ {x^3-3x \over 2}+\sum^{n-1}_{k=1} (k+1)(2k+1)(4k+3) p_{2k+1}(x)\Big], \\
\label{E413} \Gb(x) &=&  {1\over (n+1)(2n+3)}\Big[1-\sum^n_{k=1} k(2k+1)(4k+1) p_{2k}(x)\Big],\\
\nnb  a(x) &=&  - {x^3-3x \over 2}P_1(x)-\sum^{2n}_{k=2} (2k+1) p_k(x)\, P_k(x)\\
\nnb && +{1\over (n+1)(2n+1)}\Big[ {x^3-3x \over 2}+\sum^{n-1}_{k=1} (k+1)(2k+1)(4k+3) p_{2k+1}(x)\Big] P_{2n+1}(x)\\
\label{E414} && -{1\over (n+1)(2n+3)}\Big[1-\sum^n_{k=1} k(2k+1)(4k+1) p_{2k}(x)\Big]P_{2n+2}(x), \\
\nnb  {1\over 2c(x)} &=& {x^2+1\over 4}-{a(x)\over 2}\\
\nnb  &=& {(x^2-1)^2\over 4}+{1\over 2}\sum^{2n}_{k=2} (2k+1) p_k(x)\, P_k(x)\\
\nnb  && -{1\over 2(n+1)(2n+1)}\Big[ {x^3-3x \over 2}+\sum^{n-1}_{k=1} (k+1)(2k+1)(4k+3) p_{2k+1}(x)\Big]P_{2n+1}(x)\\
\label{E415} && +{1\over 2(n+1)(2n+3)}\Big[1-\sum^n_{k=1} k(2k+1)(4k+1) p_{2k}(x)\Big]P_{2n+2}(x) .
\end{eqnarray}
Finally, \refeq{E118} and \refeq{E119} follow from direct calculations. We get the proof. \endpf

\textbf{Proof of Corollary \ref{T105}.} 

By \refeq{E107} and \refeq{E109}, we have
\begin{eqnarray}\label{E416}
\nnb Q_{2n+2}(t) &=& a(x)+ {x^3-3x \over 2}P_1(t)+\sum^{2n+1}_{k=2} (2k+1) p_k(x)\, P_k(t)\\
&& +\Ga(x) P_{2n+2}(t)+\Gb(x) P_{2n+3}(t).
\end{eqnarray}
Then by  \refeq{E111},
\begin{eqnarray}
\label{E417} \Ga(x) &=& {1\over (n+1)(2n+3)}\Big[-1+\sum^n_{k=1} k(2k+1)(4k+1) p_{2k}(x)\Big],\\
\label{E418} \Gb(x) &=&-{1\over (n+2)(2n+3)}\Big[ {x^3-3x \over 2}+\sum^n_{k=1} (k+1)(2k+1)(4k+3) p_{2k+1}(x)\Big], \\
\nnb  a(x) &=&  - {x^3-3x \over 2}P_1(x)-\sum^{2n+1}_{k=2} (2k+1) p_k(x)\, P_k(x)\\
\nnb &&- {1\over (n+1)(2n+3)}\Big[-1+\sum^n_{k=1} k(2k+1)(4k+1) p_{2k}(x)\Big] P_{2n+2}(x)\\
\nnb && +{1\over (n+2)(2n+3)}\Big[ {x^3-3x \over 2}+\sum^n_{k=1} (k+1)(2k+1)(4k+3) p_{2k+1}(x)\Big]P_{2n+3}(x), \\
\label{E419} && \\
\nnb {1\over 2c(x)} &=& {x^2+1\over 4}-{a(x)\over 2}\\
\nnb  &=& {(x^2-1)^2\over 4}+{1\over 2}\sum^{2n+1}_{k=2} (2k+1) p_k(x)\, P_k(x)\\
\nnb &&+ {1\over 2(n+1)(2n+3)}\Big[-1+\sum^n_{k=1} k(2k+1)(4k+1) p_{2k}(x)\Big] P_{2n+2}(x)\\
\nnb && -{1\over 2(n+2)(2n+3)}\Big[ {x^3-3x \over 2}+\sum^n_{k=1} (k+1)(2k+1)(4k+3) p_{2k+1}(x)\Big]P_{2n+3}(x). \\
\label{E420} &&
\end{eqnarray}
Finally, \refeq{E120} and \refeq{E121} follow from direct calculations. We get the proof. \endpf

\textbf{Proof of Corollary \ref{T106}.}   First, we get
we get \refeq{E123} from  \refeq{E113}, \refeq{E115}, \refeq{E118} and \refeq{E121}.

By \refeq{E307}, $p_{2k+1}(0)=0$ ($k=0,1,2, \ldots$). Thus, if $n\geq 2$,  we get from \refeq{E117} and \refeq{E120} that
\begin{eqnarray}\label{E421}
\nnb  && B_{2n+1}(0)  =B_{2n+2}(0)\\
 \nnb &=& \Big\{ {1\over 4}+{1\over 2}\sum^{n}_{k=1} (4k+1) p_{2k}(0)\, P_{2k}(0)\\
  && +{1\over 2(n+1)(2n+3)}\Big[1-\sum^n_{k=1} k(2k+1)(4k+1) p_{2k}(0)\Big]P_{2n+2}(0) \Big\}^{1\over 2}.
\end{eqnarray}
Moreover, using \refeq{E306}--\refeq{E307}, we get \refeq{E124}:
\begin{eqnarray*}
\nnb  && B_{2n+1}(0)  =B_{2n+2}(0)\\
\nnb &=& \Big\{ {1\over 4}+{1\over 2} \times {5\over 8}\times (-{1\over 2})-{1\over 2}\sum^{n}_{k=2}(4k+1) {(2k-3)!!\over (2k+2)!!} {(2k-1)!!\over (2k)!!}\\
\nnb && -{(-1)^n\over 2(n+1)(2n+3)}\Big[1+\sum^n_{k=1} (-1)^k k(2k+1)(4k+1) {(2k-3)!!\over (2k+2)!!}\Big]{(2n+1)!!\over (2n+2)!!} \Big\}^{1\over 2}\\
\nnb &=& \Big\{ {3\over 32} -{1\over 2}\sum^n_{k=2} (4k+1) {(2k-3)!!\over (2k+2)!!}\cd{(2k-1)!!\over (2k)!!}+{7\cd (-1)^n\over 16(n+1)(2n+3)}{(2n+1)!!\over (2n+2)!!}\\
\nnb &&-\sum^n_{k=2}{(-1)^{n-k}k(2k+1)(4k+1)\over 2(n+1)(2n+3)}{(2k-3)!!\over (2k+2)!!}\cd{(2n+1)!!\over (2n+2)!!}\Big\}^{1\over 2},\q n\geq 2.
\end{eqnarray*}
And \refeq{E125} follows directly from \refeq{E124}.\endpf


\begin{thebibliography}{99}
\bibitem{Ber} L. D. Berkovitz, Optimal Control Theory, Springer-Verlag, New York, 1983.
\bibitem{Kal} G. A. Kalyabin,  Sharp constants in one-dimensional inequalities of Poincar\'e type (Russian),  Mat., 90(2011), no. 4, pp. 634--636;  translation in  Math. Notes, 90(2011),  no. 3--4, pp. 615--618.
\bibitem{Yo} K. Yosida, Functional Analysis, Springer-Verlag, 6th edition, 1980.
\end{thebibliography}
\end{document}